\documentclass[smallcondensed,referee]{svjour3}
\usepackage{lineno,fullpage}
\usepackage{amsmath,amssymb}

\usepackage{graphicx}

\newcommand*\patchAmsMathEnvironmentForLineno[1]{
 \expandafter\let\csname old#1\expandafter\endcsname\csname #1\endcsname
 \expandafter\let\csname oldend#1\expandafter\endcsname\csname end#1\endcsname
 \renewenvironment{#1}
 {\linenomath\csname old#1\endcsname}
 {\csname oldend#1\endcsname\endlinenomath}}
\newcommand*\patchBothAmsMathEnvironmentsForLineno[1]{
 \patchAmsMathEnvironmentForLineno{#1}
 \patchAmsMathEnvironmentForLineno{#1*}}
\AtBeginDocument{
\patchBothAmsMathEnvironmentsForLineno{equation}
\patchBothAmsMathEnvironmentsForLineno{align}}

\journalname{Numerical Algorithms}

\begin{document}

\title{Review of Inverse Laplace Transform Algorithms for
  Laplace-Space Numerical Approaches} \titlerunning{Numerical Inverse
  Laplace Transform Review}

\author{Kristopher L. Kuhlman}

\institute{Repository Performance Department, Sandia National Laboratories\\
  4100 National Parks Highway, Carlsbad, New Mexico, USA\\
  Tel: +1 575-234-0084\\
  \email{klkuhlm@sandia.gov}
}

\maketitle

\begin{abstract}
  A boundary element method (BEM) simulation is used to compare the
  efficiency of numerical inverse Laplace transform strategies,
  considering general requirements of Laplace-space numerical
  approaches.  The two-dimensional BEM solution is used to solve the
  Laplace-transformed diffusion equation, producing a time-domain
  solution after a numerical Laplace transform inversion.  Motivated
  by the needs of numerical methods posed in Laplace-transformed
  space, we compare five inverse Laplace transform algorithms and
  discuss implementation techniques to minimize the number of
  Laplace-space function evaluations. We investigate the ability to
  calculate a sequence of time domain values using the fewest
  Laplace-space model evaluations.  We find Fourier-series based
  inversion algorithms work for common time behaviors, are the most
  robust with respect to free parameters, and allow for
  straightforward image function evaluation re-use across at least a
  log cycle of time.


\end{abstract}

\begin{linenumbers}

\section{Introduction}
Simulation methods that are posed in Laplace-transformed space, then
numerically inverted back to the time domain (i.e., Laplace-space
methods), are a viable alternative to the more standard use of finite
differences in time.  We use the the two-dimensional boundary element
method (BEM) as an example of this type of approach, to solve the
Laplace-transformed diffusion equation (i.e., the Yukawa or modified
Helmholtz equation). We investigate five numerical inverse Laplace
transform methods and implementation approaches, namely the methods of
\cite{stehfest70}, \cite{schapery62}, \cite{weeks1966numerical},
\cite{talbot1979accurate}, and \cite{deHoog82}.  Naively implemented
Laplace-space simulations can be more computationally expensive than
using finite differences in time, but they have the advantage of
allowing evaluation at any time, without evolving from an initial
condition, and image function calculations are trivially parallelized
across Laplace parameters \cite{davies02}.  When Laplace-space
numerical models are used in parameter estimation, hundreds or
thousands of forward simulations may be required -- making forward
model efficiency critical.  Although parameter estimation may be done
directly in Laplace space \cite{barnhart12}, choosing an efficient
inversion strategy is important in most applications.

The Laplace transform has a long history of use to derive analytical
solutions to diffusion and wave problems (e.g., see list of citations
by \cite[pp.\ 191-220]{duffy2004transform}).  Often the analytical
inverse transform is too difficult to find or evaluate in closed form.
A researcher then resorts to approximate analytical methods (e.g.,
\cite{hantush1960modification,sternberg1969flow}) or numerical
inversion (e.g.,
\cite{malama09unconfinedSP,mishra2010improved}). Numerical methods can
similarly benefit from the Laplace transform, converting the
time-dependence of a differential equation to parameter dependence.
Laplace-space finite-element approaches have seen application to
groundwater flow and solute transport (e.g.,
\cite{sudicky1992laplace,morales2006laplace}), and Laplace-space BEM
has also been used in groundwater applications (e.g., \cite[\S
10.3]{kythe95} or \cite[\S 10.1]{liggett82}). The Laplace transform
analytic element method \cite{kuhlman2009laplace} is a transient
extension of the analytic element method.  These different
Laplace-space approaches may have diverse spatial solution strategies,
but they have a common requirement of effective Laplace transform
numerical inversion algorithms.  We couple a BEM model in the Laplace
domain with a numerical Laplace transform inversion routine, but our
conclusions should be valid for both gridded and mesh-free
Laplace-space numerical methods.  Any Laplace-space numerical approach
begins with determination of optimal Laplace parameter values. Then
each image function evaluation is computed from the simulation.  The
final step involves approximating the time-domain solution from the
vector of image function values using the algorithm of choice.

\cite{bellman1966} was an early review book on numerical Laplace
transform inversion for linear and non-linear problems, but without
the benefit of the many algorithms that have since been
developed. \cite{davies79} performed a thorough survey, assessing
numerical Laplace transform inversion algorithm accuracy for
techniques available in 1979, using simple functions for their
benchmarks.  \cite{duffy93} reviewed the numerical inversion
characteristics for more pathological time behaviors using the Fourier
series, Talbot, and Weeks inversion methods.  The review book by
\cite{cohen07} summarizes historical reviews and discusses commonly
used inversion their variations.  More details and examples can be
found in these reference regarding the convergence behavior of the
five inversion algorithms discussed here.

While these published numerical inverse Laplace transform algorithm
reviews are thorough and useful, they focus on computing a single
time-domain solution as accurately as possible.  These reviews did not
try to minimize Laplace-space function evaluations, since their
functions were simple closed-form expressions, not simulations. We
investigate Laplace transform inversions algorithms that can compute a
sequence of time domain values using the fewest Laplace-space model
evaluations possible, a desirable property for use in Laplace-space
numerical methods.  Using numerical Laplace transform inversion in a
simulation approach, rather than a time-marching method, allows the
researcher to readily switch between fast and accurate by changing the
number of approximation terms in the inversion.

In the next section we define the mathematical formulation of the
governing equation and Laplace transform.  In the third section we
introduce the five inverse Laplace transform algorithms.  In the final
section we compare results using five different inversion algorithms
to invert the BEM modified Helmholtz solution on the same domain with
four different boundary conditions, leading to recommendations for
Laplace-space numerical approaches.

\section{Governing Equation and Laplace Transform}
The BEM model generally simulates \textit{substance} flow (e.g.,
energy or groundwater), which can be related to a potential $\phi$
(e.g., temperature or hydraulic head).  The medium property $\alpha$
is diffusivity $[\mathrm{L^2/T}]$, the ratio of the conductance in the
substance flux and potential gradient relation (e.g., Fourier's or
Darcy's law) to the substance capacity per unit mass (e.g., heat
capacity or storativity). The BEM (e.g.,
\cite{kythe95,liggett82,brebbia84}) can be used to solve the diffusion
equation
\begin{equation}\nabla^2 \phi = \frac{1}{\alpha} \frac{\partial \phi}{\partial t},
\label{eq:Diffusion}
\end{equation}
where $\alpha$ is a real constant in space and time. We
consider \eqref{eq:Diffusion} in a domain subject to a combination of
Dirichlet $ \phi \left(\Gamma_\mathrm{u}(s),t \right)=f_\mathrm{u}(s,t)$
and Neumann $\hat{n} \cdot \nabla \phi\left(\Gamma_\mathrm{q}(s),t
\right)=f_\mathrm{q}(s,t)$ boundary conditions along the perimeter of
the 2D domain $\Gamma=\Gamma_\mathrm{u}\cup\Gamma_\mathrm{q}$, where
$\hat{n}$ is the boundary unit normal, and $s$ is a boundary
arc-length parameter.  Without loss of generality, we only consider
homogeneous initial conditions.

The Laplace transform is
\begin{equation}
  \mathcal{L}\lbrace f(t)\rbrace \equiv \bar{f}(p) =
  \int_{0}^{\infty}f(t)e^{-pt}\, \mathrm{d} t,
\label{eq:LaplaceFwd}
\end{equation}
where $p$ is the generally complex-valued Laplace parameter, and the
over-bar denotes a transformed variable. The transformed diffusion
equation with zero initial conditions is the homogeneous Yukawa or
modified Helmholtz equation,
\begin{equation}
  \nabla^2 \overline{\phi} - q^2 \overline{\phi} = 0,
\label{eq:Helmholtz}
\end{equation}
where $q^2=p/ \alpha$. Equation \ref{eq:Helmholtz} arises in several
groundwater applications, including transient, leaky, and linearized
unsaturated flow \cite{bakker2011computational}.  The transformed
boundary conditions are $\bar{\phi}\left( \Gamma_\mathrm{u}(s)\right)
= f_\mathrm{u}(s)\bar{f}_\mathrm{t}(p)$ and $\hat{n} \cdot \nabla
\bar{\phi}(\Gamma_\mathrm{q}(s))=f_\mathrm{q}(s)\bar{f}_\mathrm{t}(p)$,
where the temporal and spatial behaviors have been decomposed as in
separation of variables. Arbitrary time behavior can be developed
through convolution in $t$ (Duhamel's theorem), which is
multiplication of image functions in Laplace space.  Here,
$\bar{f}_\mathrm{t}(p)$ represents the Laplace transform of the time
behavior applied to the boundary conditions.  The Laplace
transformation makes it possible to solve transient diffusion (a
parabolic equation) using the BEM, which is well-suited for
elliptical-type equations.

The inverse Laplace transform is defined as the Bromwich contour
integral,
\begin{equation}
\mathcal{L}^{-1}\left\lbrace\bar{f}(p)\right\rbrace = f(t) =
\frac{1}{2 \pi i}\int_{\sigma-i\infty}^{\sigma+i\infty}\bar{f}(p) e^{pt}\, \mathrm{d} p,
\label{eq:LaplaceInv}
\end{equation}
where the abscissa of convergence $\sigma>0$ is a real constant chosen
to put the contour to the right of all singularities in
$\bar{f}(p)$. In Laplace-space numerical approaches,
\eqref{eq:Helmholtz} is solved by a suitable numerical method,
therefore only samples of $\bar{f}(p)$ are available; this precludes
an analytical inversion. Five numerical inverse Laplace transform
algorithms are discussed in the following section.

\section{Numerical Inverse Laplace Transform Methods}
Equation \ref{eq:LaplaceInv} is an integral equation for unknown
$f(t)$ given $\bar{f}(p)$; its numerical solution is broadly split
into two categories. Methods are either based on quadrature or
functional expansion using analytically invertible basis functions.
\cite[Chap.\ 19]{davies05} relates most major classes of inverse
Laplace transform methods using a unified theoretical foundation; we
adopt a simplified form of their general notation.  The Fourier series
and Talbot methods are quadrature-based, directly approximating
\eqref{eq:LaplaceInv}.  Weeks' and Piessen's methods are $\bar{f}(p)$
expansions using complex-valued basis functions, while the
Gaver-Stehfest and Schapery methods use real-valued functions to
accomplish this.

The numerical inverse Laplace transform is in general an ill-posed
problem (e.g., \cite{al-shuaibi01}).  No single approach is optimal for all
circumstances and temporal behaviors, leading to the diversity of
viable numerical approaches in the literature (e.g., \cite{cohen07}).

\subsection{Gaver-Stehfest Method}
The Post-Widder formula \cite{widder41,al-shuaibi01} is an
approximation to \eqref{eq:LaplaceInv} that only requires $\bar{f}(p)$
for real $p$ to represent \eqref{eq:LaplaceFwd} as an asymptotic
Taylor series expansion.  The formula requires high-order analytic
image function derivatives, and is impractical for numerical
computation. Stehfest proposed a discrete version of the Post-Widder
formula using finite differences and Salzer summation
\cite{stehfest70},
\begin{equation}
  f(t,N)=\frac{\ln 2}{t}\sum_{k=1}^{N}V_k\bar{f}\left( k\frac{\ln 2}{t} \right).
\label{eq:GaverStehfest}
\end{equation}
The $V_k$ coefficients only depend on the number of expansion terms,
$N$ (which must be even), which are
\begin{equation}
  V_k=(-1)^{k + N/2}\sum^{\min(k,N/2)}_{j=\lfloor(k+1)/2 \rfloor}
  \frac{j^{\frac{N}{2}}(2j)!}{(\frac{N}{2}-j)!\,j!\,(j-1)!\,(k-j)!\,(2j-k)!}.
\label{eq:GaverStehfestV}
\end{equation}
These become very large and alternate in sign for increasing $k$. The
sum \eqref{eq:GaverStehfest} begins to suffer from cancellation for
$N \ge$ the number of decimal digits of precision (e.g., double
precision $=16$). For $\bar{f}_\mathrm{t}(p)$ that are non-oscillatory
and continuous, $N\le18$ is usually sufficient \cite{stehfest70}. If
programmed using arbitrary precision (e.g. Mathematica or a
multi-precision library \cite{bailey02,mpmath}), the method can be
made accurate for most cases \cite{abate04}.  Unfortunately, $p$ is
explicitly a function of $t$; for each new $t$, a new $\bar{f}(p)$
vector is needed. In Laplace-space numerical approaches, each vector
element is constructed using a simulation, therefore this can be a
large penalty.

The method is quite easy to program; the $V_j$ can be computed once
and saved as constants. This method has been popular due to its
simplicity and adequacy for exponentially decaying
$\bar{f}_\mathrm{t}(p)$.

\subsection{Schapery's Method}
We can expand the deviation of $f(t)$ from steady-state $f_\mathrm{s}$
using exponential basis functions \cite{schapery62},
\begin{equation}
f(t,N) = f_\mathrm{s}+\sum_{i=1}^{N}a_ie^{-p_it} ,
\label{eq:SchaperyForm}
\end{equation}
where $a_i$ is a vector of unknown constants. Applying
\eqref{eq:LaplaceFwd} to \eqref{eq:SchaperyForm} gives
\begin{equation}
  \bar{f}(p_j,N) = \frac{f_\mathrm{s}}{p_j}+\sum_{i=1}^{N}\frac{a_i}{p_i+p_j}
  \qquad j=1,2,\ldots,M.
\label{eq:SchaperyTransformed2}
\end{equation}
The $p_j$ are selected (a geometric series is recommended
 \cite{liggett82}) to cover the important fluctuations in $\bar{f}(p)$.
After setting $p_i=p_j$ the $a_i$ coefficients can be determined as
the solution to $P_{ij}a_i=\left(\bar{f}(p_j) -
  f_\mathrm{s}/p_j\right)$. The symmetric matrix to decompose is
\begin{equation*}
  P_{ij}=\left[
    \begin{matrix}
      (2p_1)^{-1} & (p_1+p_2)^{-1} & \ldots & (p_1+p_N)^{-1} \\
      (p_2+p_1)^{-1} & (2p_2)^{-1} & \ldots & (p_2+p_N)^{-1} \\
      \vdots        & \vdots        & \ddots & \vdots \\
      (p_N+p_1)^{-1} & (p_N+p_2)^{-1} & \ldots & (2p_N)^{-1}
    \end{matrix}
  \right],
  \label{eq:SchaperyMatrixP}
\end{equation*}
which only depends on $p_j$; it can be decomposed
independently of $\bar{f}(p)$ and $f_\mathrm{s}$.

This method is not difficult to implement when existing matrix
decomposition libraries are available, and only requires real
computation. The method has been used for inverting BEM results
\cite{liggett82}, but has two main drawbacks. First, in its
formulation above, it requires a steady-state solution, but
\eqref{eq:SchaperyForm} could be posed without $f_\mathrm{s}$.
Secondly, no theory is presented for optimally picking $p_j$; some
trial and error is required \cite{liggett82}.

\subsection{M\"{o}bius Transformation Methods}
We can use the M\"{o}bius transformation to conformally map the
half-plane right of $\sigma$ to the unit disc, making the Laplace
domain more amenable to approximation using orthonormal polynomials
(e.g., Chebyshev \cite{piessens72}, \cite[\S 28]{lanczos88} or
Laguerre \cite{weeks1966numerical}, \cite{lyness86}, \cite[\S
30]{lanczos88}).  If $\sigma$ was chosen properly, $\bar{f}(p)$ is
guaranteed to be analytic inside the unit circle.  The most-used
inverse Laplace transform method from this class is Weeks' method,
which uses a complex power series to expand $\bar{f}(p)$ inside the
unit circle. Upon inverse Laplace transformation, the power series
becomes a Laguerre polynomial series.

Weeks method is 
\begin{equation}
f(t,N+1) = e^{\left(\kappa-b/2\right)t}
\sum_{n=0}^{N}a_n \text{L}_n(bt),
\label{eq:Weeks}
\end{equation}
where $\text{L}_n(z)$ is an $n$-order Laguerre polynomial and $\kappa$
and $b$ are free parameters. Weeks suggested $\kappa=\sigma +
1/t_\mathrm{max}$ and $b=N/t_\mathrm{max}$, where $t_\mathrm{max}$ is
the maximum transformed time.  The parameters $b$ and $\kappa$ are
chosen to optimize convergence; some schemes are given
\cite{weideman99} for finding optimum parameter values for a given
$\bar{f}_\mathrm{t}(p)$, but search techniques require hundreds of
$\bar{f}(p)$ evaluations.  A more general form of \eqref{eq:Weeks} can
also be used, which allows for more general asymptotic behavior of the
image function \cite[\S 19.5]{davies05}.  Weeks assumed $p\bar{f}(p)$
is analytic at infinity. The Laplace transform of \eqref{eq:Weeks} is
known, but to make it easier to represent with polynomials,
$\bar{f}(p)$ is mapped inside the unit circle via
$z=(p-\kappa-2b)/(p-\kappa+2b)$. The coefficients $a_n$ are determined
by the midpoint rule,
\begin{equation}
a_n = \frac{1}{2M} \sum_{j=-M}^{M-1}
  \Psi\left[ \exp\left(i\theta_{j-\frac{1}{2}}\right) \right]
  \exp\left(-in\theta_{j-\frac{1}{2}}\right)
\label{eq:WeeksCoeff}
\end{equation}
where $\theta_j=j\pi/M$ and the conformally-mapped image
function is
\begin{equation}
\Psi(z)=\frac{b}{1-z}\bar{f}\left( \kappa -\frac{b}{2} + \frac{b}{1-z}  \right).
\label{eq:WeeksMappedFp}
\end{equation}
The argument of $\bar{f}(z)$ in \eqref{eq:WeeksMappedFp} is the
inverse mapping of $z \mapsto p$, it shows $p$ does not functionally
depend on $t$, but Weeks' rules-of-thumb for $b$ and $\kappa$ depend on
$t_\mathrm{max}$.

There are other related methods which use different orthonormal
polynomials to represent $\bar{f}(p)$ inside the unit circle.
Chebyshev polynomials (known as Piessen's method \cite{piessens72})
can be used to expand the $\bar{f}(z)$ on the real interval
$[-1,1]$. The Weeks method is moderately easy to program, requiring
the use of Clenshaw recurrence formula to accurately implement
Laguerre polynomials.  Piessen's method is similar to implement, with
a similar recurrence formula for Chebyshev polynomials.

\subsection{Talbot Method}
We can deform the Bromwich contour into a parabola around the negative
real axis if $\bar{f}(p)$ is analytic in the region between the
Bromwich and the deformed Talbot contours
 \cite{talbot1979accurate}. Numerically, $\bar{f}(p)$ must not overflow
as $p \rightarrow -\infty$ (e.g., in the BEM implementation, the
Green's function is the second-kind modified Bessel function, which
grows exponentially as $p\rightarrow -\infty$).  Oscillatory
$\bar{f}_\mathrm{t}(p)$ often have pairs of poles near the imaginary
$p$ axis; these poles must remain to the left of the deformed contour.

The Talbot method makes the Bromwich contour integral converge
rapidly, since $p$ becomes large and negative, making the $e^{pt}$
term in \eqref{eq:LaplaceInv} very small. A one-parameter ``fixed''
Talbot method was implemented \cite{abate04}; the Bromwich contour is
parametrized as $p(\theta)=r\theta(\cot(\theta) +i)$, where
$0\le\theta\le\pi$, and as a rule of thumb
$r=2M/(5t_\mathrm{max})$. The fixed Talbot method is
\begin{equation}
f(t,N) = \frac{r}{N}\left[ \frac{\bar{f}(r)}{2}e^{rt} + \sum_{k=1}^{N-1}
  \Re\left\lbrace e^{tp(\theta_k)}\bar{f}\left[p(\theta_k)\right]
    \left[1+i \zeta(\theta_k) \right]\right\rbrace\right],
\label{eq:fixedTalbot}
\end{equation}
where $\zeta(\theta_k) = \theta_k +
\left[\theta_k\cot(\theta_k)-1\right]\cot(\theta_k)$ and
$\theta_k=k\pi/N$ \cite{abate04}.  Although $\bar{f}(p)$
doesn't depend on $t$, the free parameter $r$ depends on
$t_\mathrm{max}$.

Step change $\bar{f}_\mathrm{t}(p)$ for non-zero time become very
large as $p \rightarrow -\infty$, since
$\mathcal{L}\left[H(t-\tau)\right]=e^{-\tau p}/p$, where $H(t-\tau)$
is the Heaviside step function centered on time $\tau$. This can lead
to precision loss, and stability or convergence issues with the
underlying numerical model, although Mathematica's arbitrary precision
capabilities have been used to get around this problem \cite{abate04}.

The fixed Talbot method is very simple to program; \cite{abate04}
provide a ten-line Mathematica implementation.

\subsection{Fourier Series Method}
We can manipulate \eqref{eq:LaplaceInv} into a Fourier
transform; first it is expanded into real and imaginary parts ($
p=\gamma+i\omega $),
\begin{equation*}
  f(t)=\frac{e^{\gamma t}}{2\pi i}\int_{-\infty}^{\infty} \left[\cos(\omega t) +
  i\sin(\omega t)\right] \left\lbrace\Re\left[\bar{f}(p)\right] +
  i\Im\left[\bar{f}(p)\right] \right\rbrace i \, \mathrm{d}\omega.
\label{eq:BromwichFourierWide}
\end{equation*}
Multiplying out the terms, keeping only the real part due
to $f(t)$ symmetry, and halving the integration range due to
symmetry again, leaves
\begin{equation}
  f(t)=\frac{e^{\gamma t}}{\pi}\int_{0}^{\infty}
  \Re\left[\bar{f}(p)\right]\cos(\omega t) -
  \Im\left[\bar{f}(p)\right]\sin(\omega t)\, \mathrm{d}\omega.
\label{eq:BromwichFourierReal}
\end{equation}
When $f(t)$ is real, \eqref{eq:BromwichFourierReal} can be represented
using the complex form or just its real or imaginary parts.  Although
these three representations are equivalent, when evaluating
\eqref{eq:BromwichFourierReal} with the trapezoid rule, the full
complex form gives the smallest discretization error
\cite{davies05}. The trapezoid rule approximation to
\eqref{eq:BromwichFourierReal} is essentially a discrete Fourier
transform,
\begin{equation}
  f(t,2N+1) = \frac{e^{\gamma t}}{T} \sum_{k=0}^{2N}{}^{'}
  \Re\left[\bar{f}\left(\gamma_0 +\frac{i\pi k}{T}\right)
    \exp\left(\frac{i\pi t}{T}\right)\right],
\label{eq:BromwichFourierTrap}
\end{equation}
where $\gamma_0=\sigma - \log (\epsilon)/T$, $\epsilon$ is the
desired relative accuracy (typically $10^{-8}$ in double precision),
$T$ is a scaling parameter (often $2t_\mathrm{max}$), and the prime
indicates the $k=0$ summation term is halved.  The $p$ in
\eqref{eq:BromwichFourierTrap} do not depend on $t$, but the free
parameter $T$ depends on $t_\mathrm{max}$.

The non-accelerated Fourier series inverse algorithm
\ref{eq:BromwichFourierTrap} is almost useless because it requires
thousands of $\bar{f}(p)$ evaluations \cite[\S
9.8]{antia02}. Practical approaches accelerate the convergence of the
sum in \ref{eq:BromwichFourierTrap}.  Although this is sometimes
called a fast-Fourier transform (FFT) method (e.g., \cite[Chap
4.4]{cohen07}), rarely do the number of $\bar{f}(p)$ evaluations in an
accelerated approach justify an FFT approach.  The method implemented
uses non-linear double acceleration with Pad\'{e} approximation and an
analytic expression for the remainder in the series
\cite{deHoog82}. Although there are several other ways to accelerate
the Fourier series approach \cite{cohen07}, this method is popular and
straightforward.  Non-linear acceleration techniques drastically
reduce the required number of function evaluations, but can lead to
numerical dispersion \cite{kano05,morales2006laplace}.  For diffusion,
dispersion associated with non-linear acceleration is not noticeable.
Schapery's, Talbot's, and Weeks' methods are not accelerated in a
non-linear manner, and therefore may lead to less numerical
dispersion, which may be more important in wave systems with sharp
fronts.

The creation of the Pad\'{e} approximation \cite{deHoog82} is
relatively straightforward in programming languages that facilitate
matrix manipulations (e.g., modern Fortran, Matlab, or NumPy
 \cite{oliphant07}).  There is no dependence on matrix decomposition
routines.

\subsection{Algorithm Properties Summary}
Table~\ref{tab:analyticSummary} summarizes aspects of the five inverse
methods.  The third column indicates whether $p$ is explicitly a
function of $t$, the fourth column indicates if the rules-of-thumb
used for the optimum parameters depend on $t_\mathrm{max}$, and the
fifth column indicates whether the transform requires complex $p$ and
$\bar{f}(p)$.

\begin{table} 
  \caption{Algorithmic Summary}
  \centering
  \begin{tabular}{|c|c|c|c|c|}
    \hline
    \textbf{Method} & \textbf{Limitations on $ \bar{f}(p)$ and $ f(t)$} &
    $p(t)$? & $p(t_\mathrm{max})$?& $ p $ \\
    \hline
    Stehfest &  no oscillations, no discontinuities in $f(t)$ & yes & no & real \\
    Schapery & smoothly varying $f(t)$, $f_\mathrm{s}$ exists & no & no & real \\
    Weeks & none & no & yes & complex  \\
    fixed Talbot & no high-frequency $f(t)$, $\bar{f}(p\rightarrow -\infty)$ exists&
    no & yes & complex \\
    Fourier series & none  & no & yes & complex \\
    \hline
  \end{tabular}
  \label{tab:analyticSummary}
\end{table}

For all methods considered here, computational effort to compute
$f(t)$ from the vector of $\bar{f}(p)$ values was insignificant
compared to the effort required to compute the BEM solution used to fill
the $\bar{f}(p)$ vector.  This suggests a more complicated
method, which allows re-use of $\bar{f}(p)$ across more values of $t$
and converges in less evaluations of $\bar{f}(p)$, would be efficient
for Laplace-domain numerical methods.  If existing libraries or
simulations only support real arguments, then the
Stehfest, Schapery, or Piessen's methods must be used.  Complex $p$
methods will pay a slight penalty in computational overhead compared
to real-only $p$ routines.  Computing with arbitrary or
higher-than-double precision (e.g., \cite{abate04}) will incur a much
larger penalty than the change from real to complex double precision.
Generally, complex $p$ methods have better convergence properties than
real-only methods. Expansion of $\bar{f}(p)$ along the real $p$ axis
is separation of non-orthogonal exponentials, while expansion along
the imaginary $p$ axis is separation of oscillatory functions \cite[\S
29]{lanczos88}.

\section{Numerical Comparison}
Four test problems were solved using the BEM for values of $p$
required by each algorithm's rules of thumb.  The test problem domain
is a $3 \times 2$ rectangle, with homogeneous initial conditions and
specified potential at two ends $\bar{\phi}(x=0)=-2
\bar{f}_\mathrm{t}(p)$, and $\bar{\phi}(x=3)=2 \bar{f}_\mathrm{t}(p)$,
and zero normal flux along the other sides $\partial
\bar{\phi}/\partial y (y=\lbrace0,2\rbrace)=0$.  All plots show the
solution computed at a point closer to the $x=0$ boundary $(x=1/3)$,
midway between the insulated boundaries $(y=1)$.

The first problem computes $\bar{f}(p)$ using the optimum $p$ at each
$t$ (like most inverse Laplace transform surveys), according to the
rules-of-thumb for each method.  While this is most accurate, it is
very inefficient -- especially when many values of $t$ are
required. In the following sections, all methods except Stehfest use
the same $\bar{f}(p)$ to invert all $t$.  A method's sensitivity to
non-optimal free parameters is important in practical use for
Laplace-space numerical approaches. By inverting more than one time
with the same set of Laplace-space function evaluations, large gains
in efficiency can be made.  The $t$ range used in the plots spans
three orders of magnitude; it was chosen to show the evolution of
potential and substance flux from initial conditions to steady
state.

\subsection{Steady Boundary Conditions, Optimum $p$}
\begin{figure} 
  \includegraphics[width=0.5\textwidth]{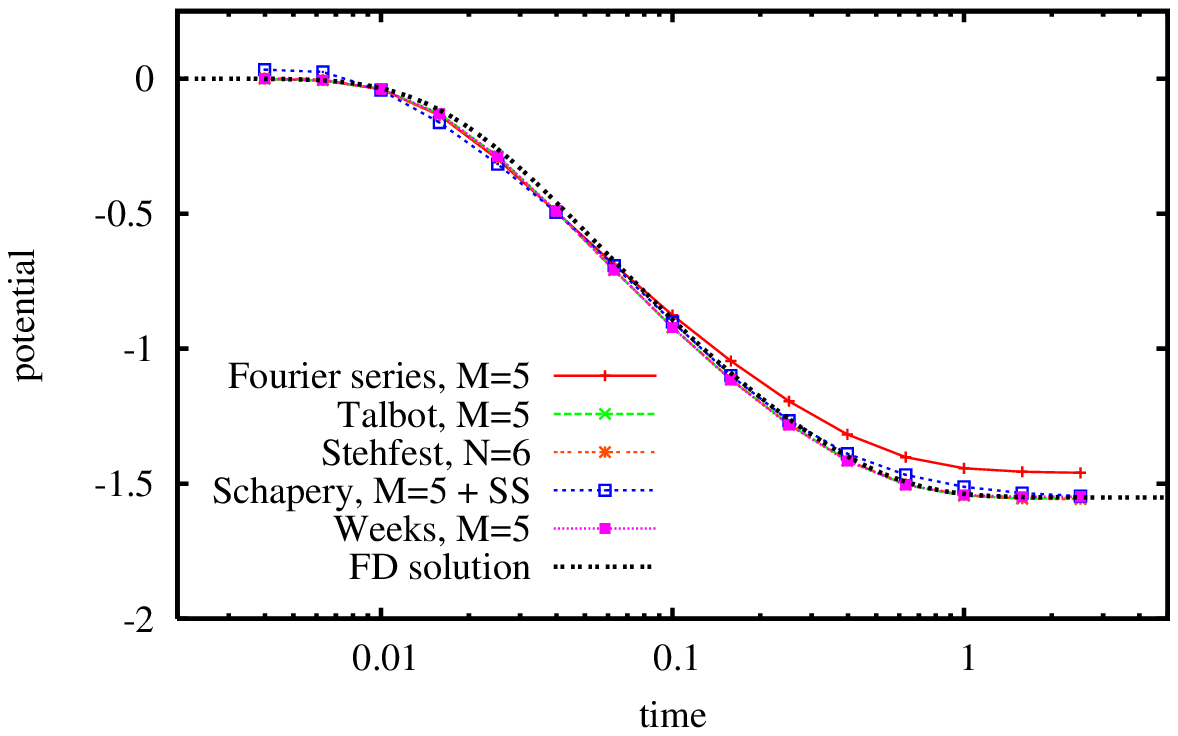}
  \includegraphics[width=0.5\textwidth]{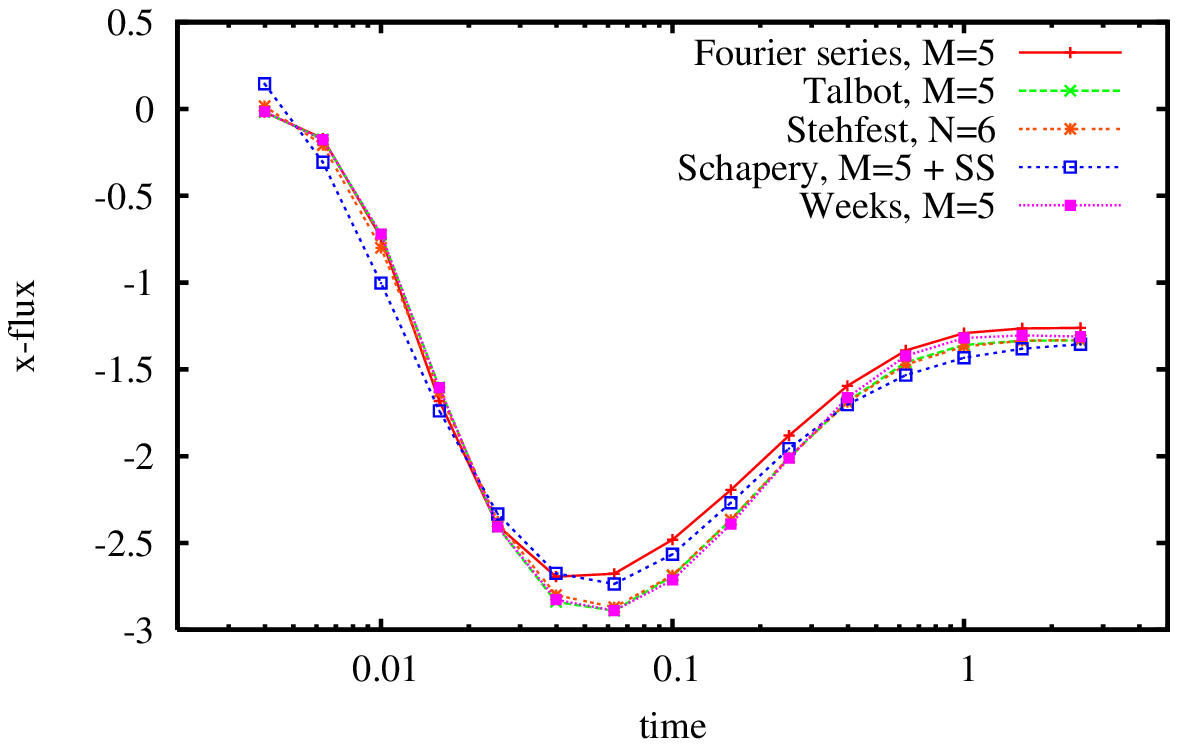}
  \caption{Plots of potential and flux through time with five methods
    for $\bar{f}_\mathrm{t}(p)=1/p$, using optimum $p$ at each
    $t$. $15 \times 5=75$ total $\bar{f}(p)$ evaluations are used by each
    method.}
  \label{fig:low-ss}
\end{figure}

The first problem has steady-state boundary conditions.  The transient
behavior is solely due to evolution from the zero initial condition,
$\bar{f}_\mathrm{t}(p)=\mathcal{L}[H(t)]=1/p$; $\bar{f}(p)$ has a
pole at the origin.  All methods performed equally well with this
simple test problem, although the Fourier series method deviates from
the finite difference solutions at larger
time. Figure~\ref{fig:low-ss} shows the inverted potential and flux
using as few evaluations of $\bar{f}(p)$ possible, without
major deviations from the finite difference benchmark solution. Some
trial and error was needed to use the Schapery method (i.e., further
optimization may be possible).

\begin{figure} 
  \includegraphics[width=0.5\textwidth]{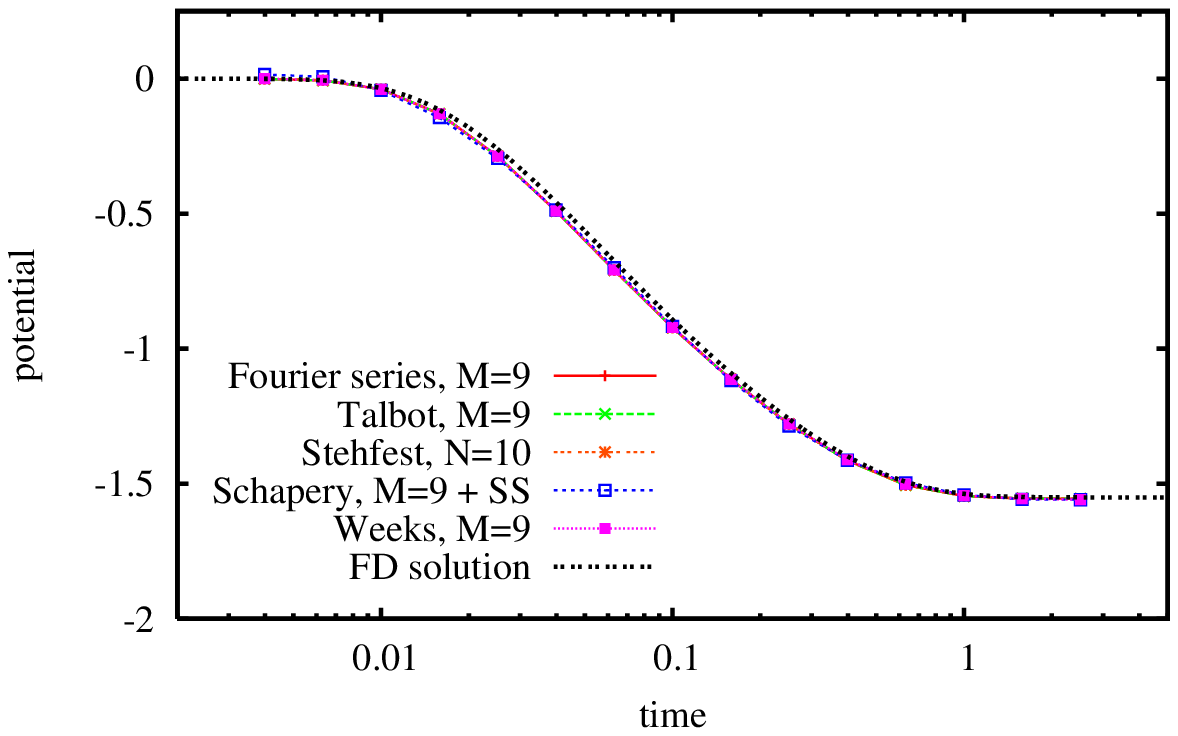}
  \includegraphics[width=0.5\textwidth]{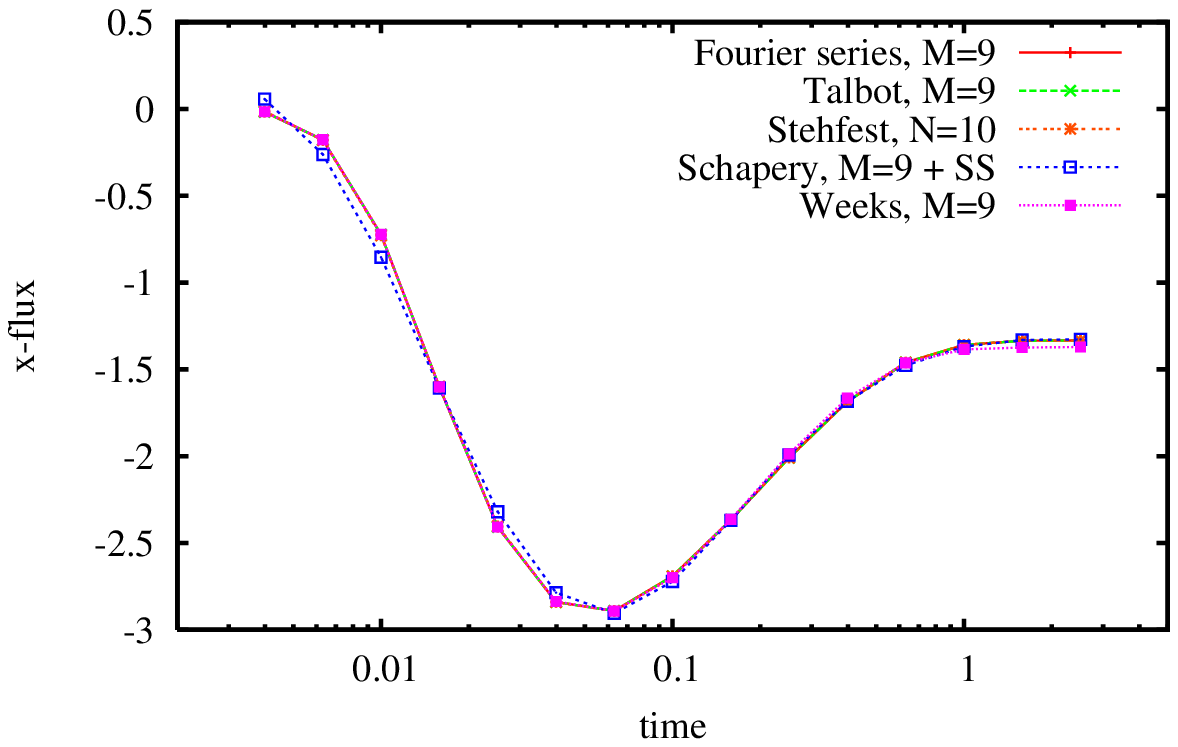}
  \caption{Plots of potential and flux through time with five methods
    for $\bar{f}_\mathrm{t}(p)=1/p$, using optimum $p$ at each
    $t$. $15 \times 9=135$ total $\bar{f}(p)$ evaluations are used by each
    method. Fourier series, Talbot, Stehfest, and Weeks curves are
    nearly coincident.}
  \label{fig:better-ss}
\end{figure}

As shown in Figure~\ref{fig:better-ss}, all the methods performed very
well for nine $\bar{f}(p)$ terms per $t$ but at least 135 $\bar{f}(p)$
evaluations are needed total for each method.  Schapery's method does
the worst in this case, but this may be improved with further
optimization of $p_j$ values.  The finite-difference approach took at
least an order of magnitude less computational effort for the given
accuracy. Making Laplace-space numerical methods useful alternatives
to traditional time-marching approaches, requires improvements to this
inefficiency.

\subsection{Steady Boundary Conditions, Same $p$}
\begin{figure} 
  \includegraphics[width=0.5\textwidth]{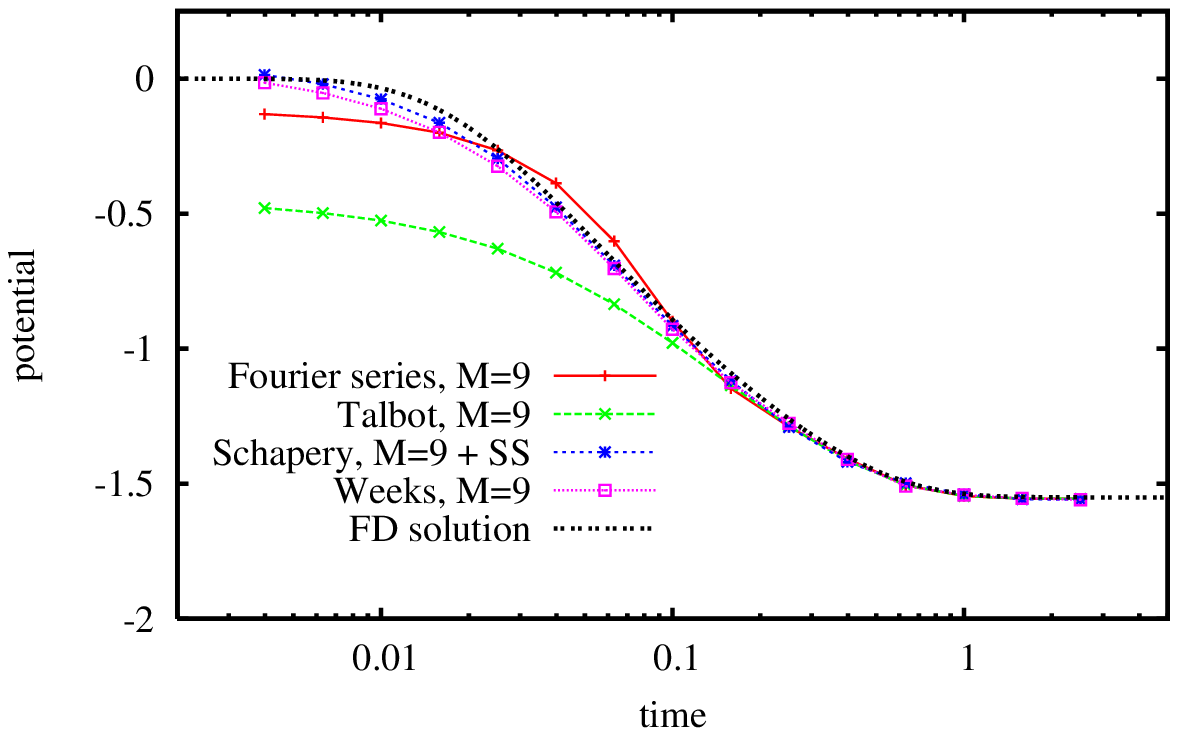}
  \includegraphics[width=0.5\textwidth]{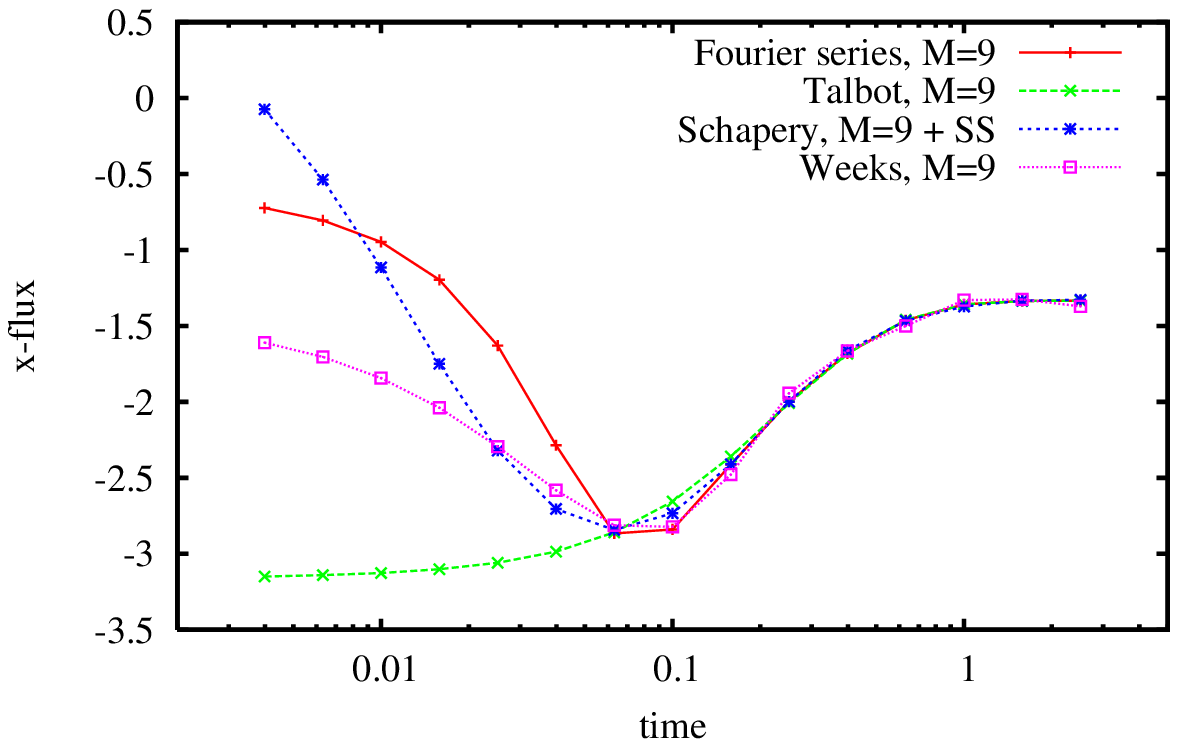}
  \caption{Plots of potential and flux through time with four methods
    for $\bar{f}_\mathrm{t}(p)=1/p$, same $p$ used across all
    $t$. Nine total $\bar{f}(p)$ evaluations are used by each method.}
  \label{fig:low-ss-rangep}
\end{figure}

All methods had more difficulty obtaining accurate results for a wide
$t$ range using only one vector of $\bar{f}(p)$ (no Stehfest
method, since $p$ explicitly depends on $t$). Only the last log-cycle
of times is inverted accurately when using nine $\bar{f}(p)$
(Figure~\ref{fig:low-ss-rangep}). All the methods -- except possibly
Schapery's -- have a more difficult time with the flux at early time
(especially the fixed Talbot method).  The apparent success of
Schapery's method can be attributed to the expansion of the deviation
from steady-state, which in this case decays exponentially with time.

\begin{figure} 
  \includegraphics[width=0.5\textwidth]{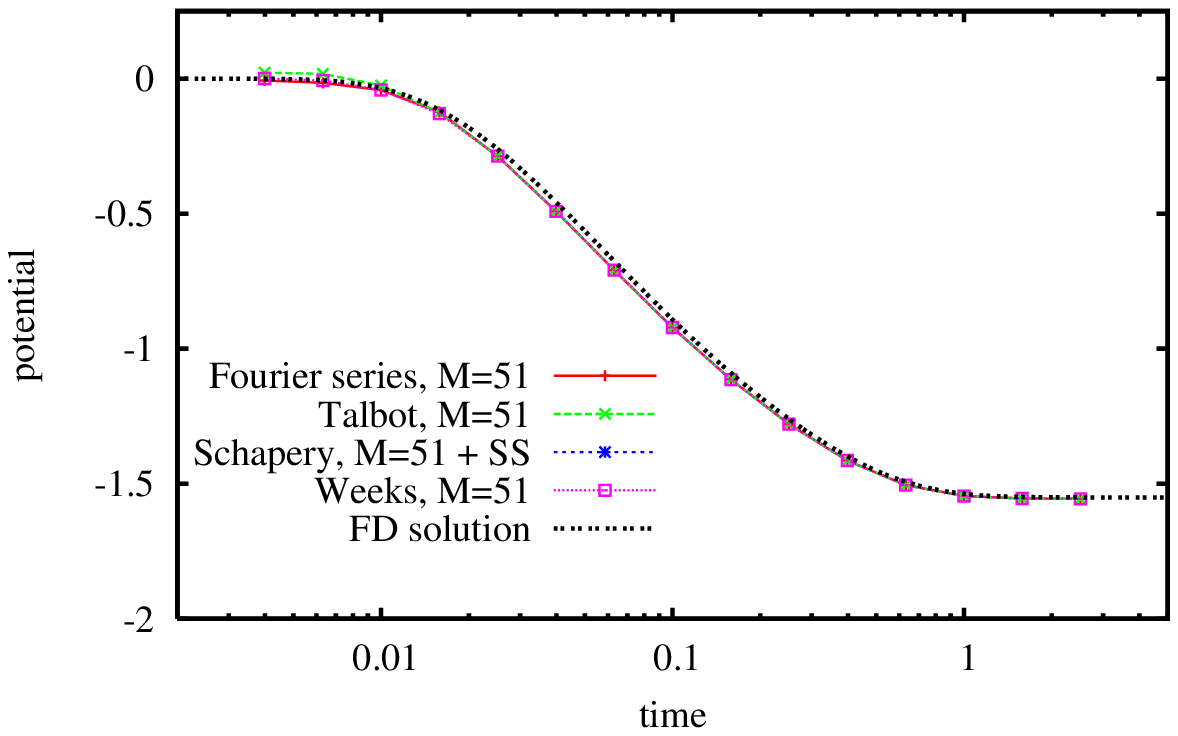}
  \includegraphics[width=0.5\textwidth]{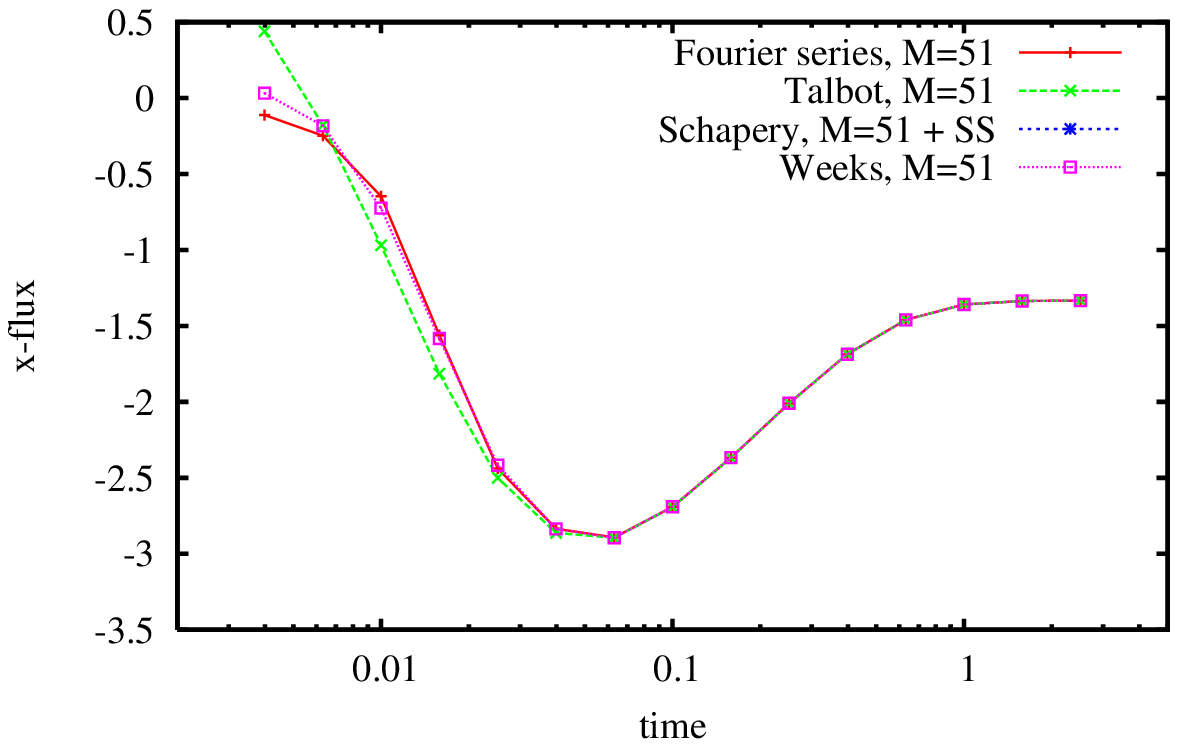}
  \caption{Plots of potential and flux through time with four methods
    for $\bar{f}_\mathrm{t}(p)=1/p$, same $p$ used across all $t$. 51
    total $\bar{f}(p)$ evaluations used by each method.  Schapery and
    Weeks curves are nearly coincident.}
  \label{fig:high-ss-rangep}
\end{figure}

Figure~\ref{fig:high-ss-rangep} shows that when increasing to 51
$\bar{f}(p)$ terms, most convergence problems disappear, except at
small times.  Grouping $t$ values by log-cycles and inverting them
together using the same $\bar{f}(p)$ is more economical than using the
optimal $p$ for each $t$ and is still relatively accurate.  The
results shown in Figure~\ref{fig:high-ss-rangep} are nearly as
accurate as those shown in Figure~\ref{fig:better-ss}, but required
1/3 the $\bar{f}(p)$ model evaluations.

\subsection{Sinusoidal Boundary Conditions, Same $p$}
This problem uses temporally sinusoidal boundary conditions,
$\bar{f}_\mathrm{t}(p)=\mathcal{L}(\cos 4t)=\frac{p}{p^2+16}$. This
boundary condition violates some assumptions of the inverse transform
algorithms (i.e., no steady-state solution, oscillatory in time), but
the behavior is still relatively simple and smooth, with singularities
at $p=\pm 4i$.

\begin{figure} 
  \includegraphics[width=0.5\textwidth]{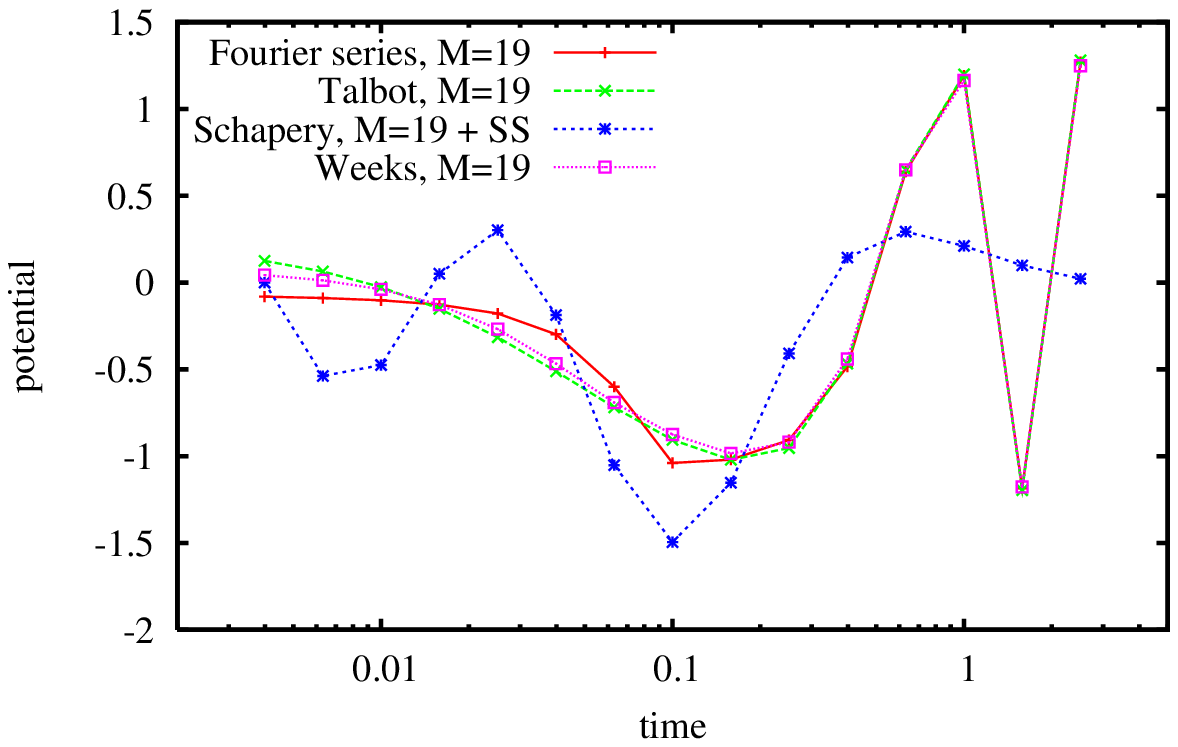}
  \includegraphics[width=0.5\textwidth]{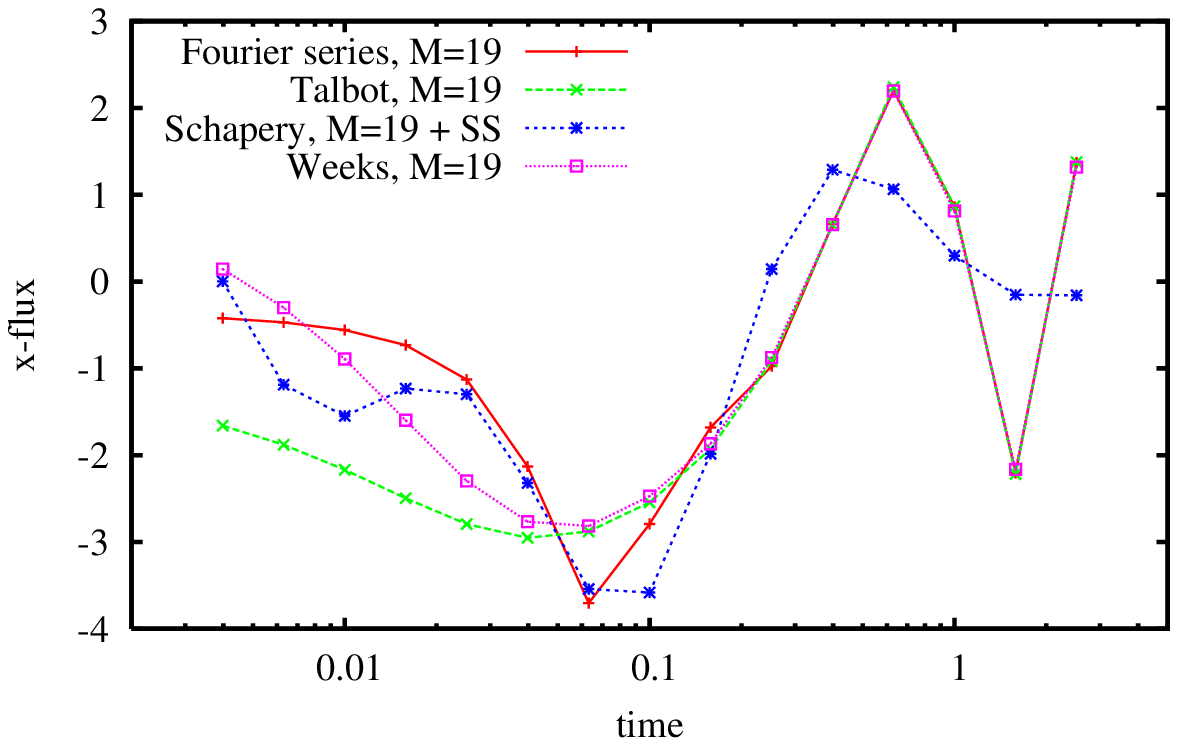}
  \caption{Plots of potential and flux through time with four methods
    for $\bar{f}_\mathrm{t}(p)=p/(p^2+16)$, same $p$ used across all
    $t$. 19 total $\bar{f}(p)$ evaluations are used by each method.}
  \label{fig:cos-low}
\end{figure}

\begin{figure} 
  \includegraphics[width=0.5\textwidth]{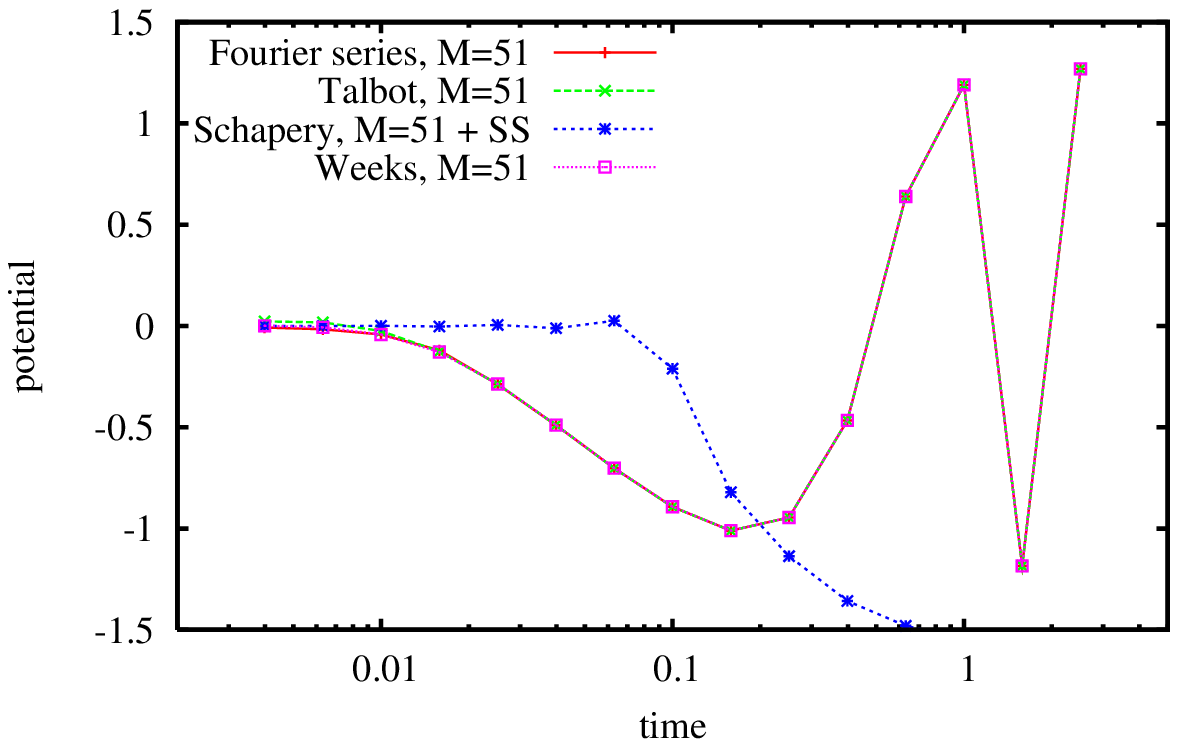}
  \includegraphics[width=0.5\textwidth]{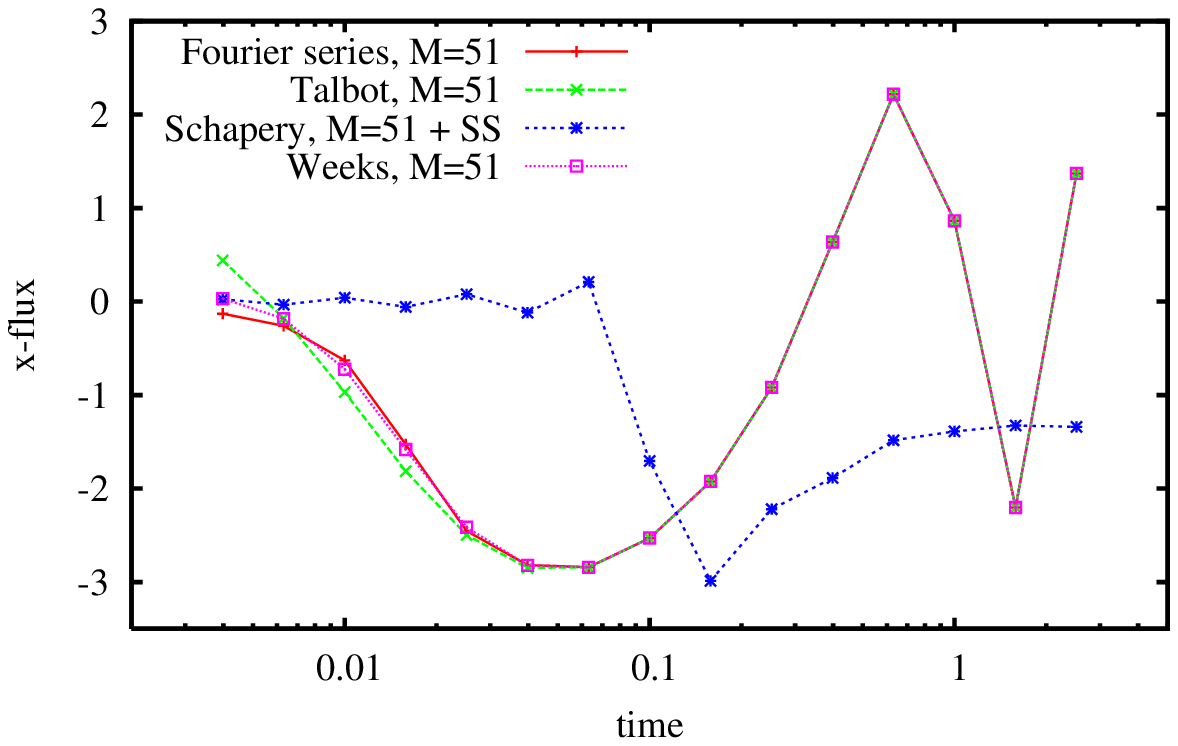}
  \caption{Plots of potential and flux through time with four methods
    for $\bar{f}_\mathrm{t}(p)=p/(p^2+16)$, same $p$ used across all
    $t$. 51 total $\bar{f}(p)$ evaluations are used by each
    method. Fourier series, Talbot, and Weeks curves are nearly
    coincident for potential.}
  \label{fig:cos-better}
\end{figure}

Figure~\ref{fig:cos-low} shows the Schapery method fails since there
is no $f_\mathrm{s}$, but the other methods do well for 19 terms
across one $t$ log cycle.  Figure~\ref{fig:cos-better} shows all
methods besides Schapery do well for 51 terms, across at least two $t$
log cycles.  A modified version of \eqref{eq:SchaperyTransformed2}
substituting $\frac{p_j}{p_j^2+16}$ for $f_\mathrm{s}/p_j$ could
extend Schapery's approach to this case, but this solution was not
considered here because of its problem specificity.

\subsection{Step-Change Boundary Condition for $\tau>0$, same $p$}
\begin{figure} 
  \includegraphics[width=0.5\textwidth]{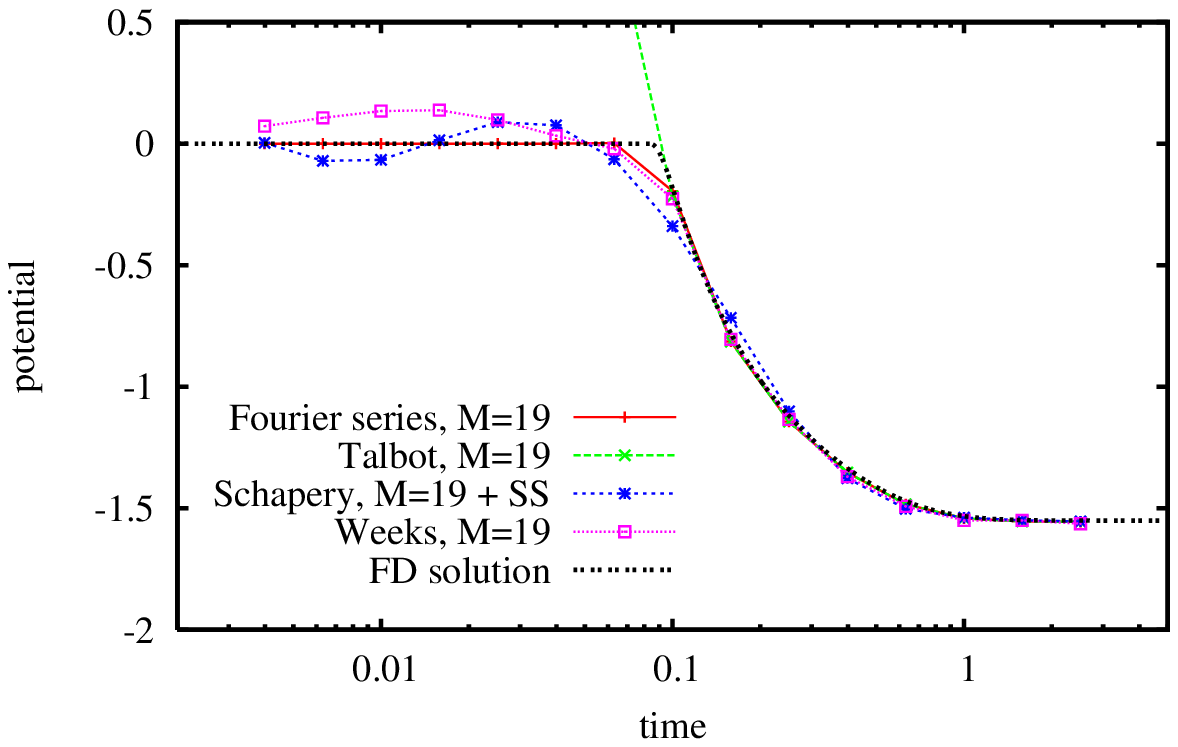}
  \includegraphics[width=0.5\textwidth]{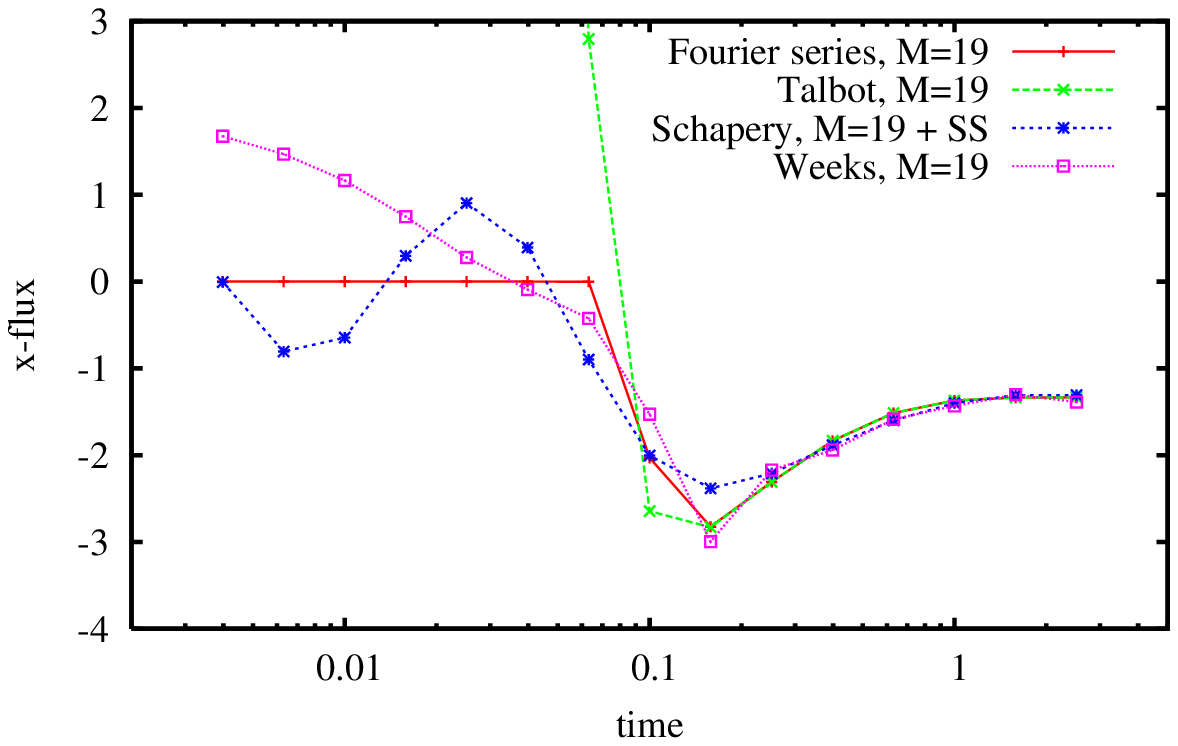}
  \caption{Plots of potential and flux through time with four methods
    for $\bar{f}_\mathrm{t}(p)=\exp(-0.08 p)/p$, same $p$
    used across all $t$. 19 total $\bar{f}(p)$ evaluations are used by each
    method.  Weeks' solution is undefined for $t<0.08$.}
  \label{fig:step-low}
\end{figure}

\begin{figure} 
  \includegraphics[width=0.5\textwidth]{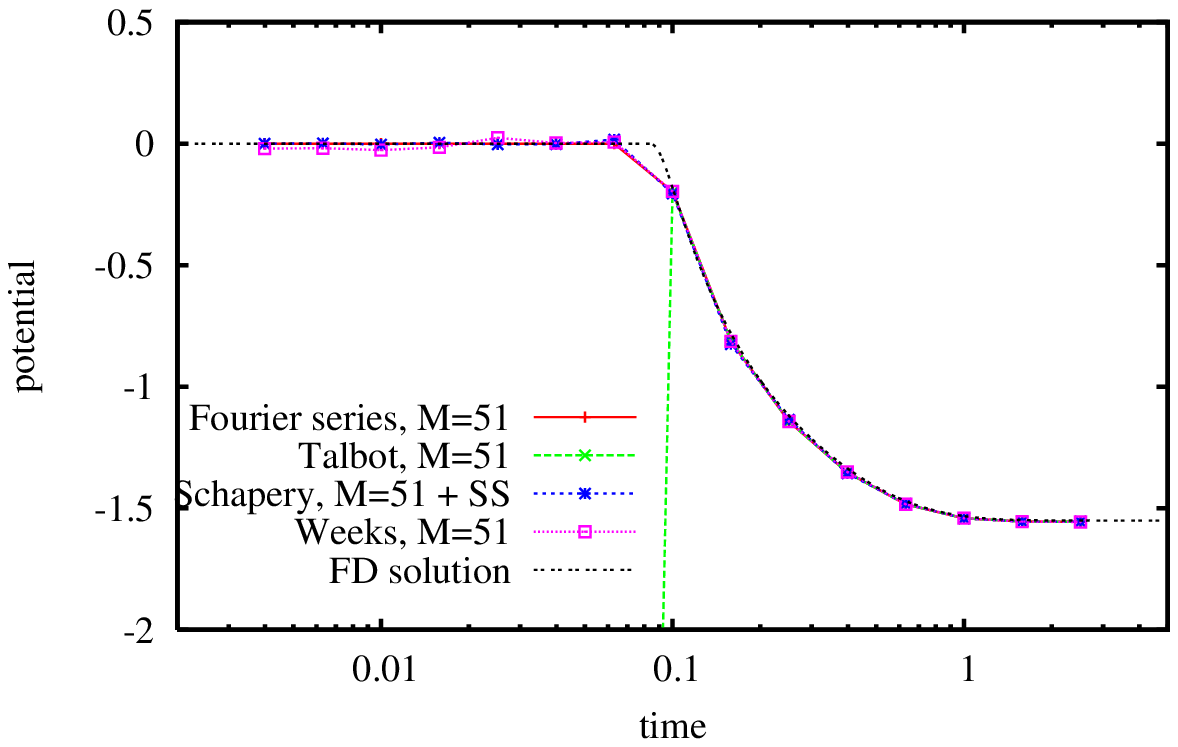}
  \includegraphics[width=0.5\textwidth]{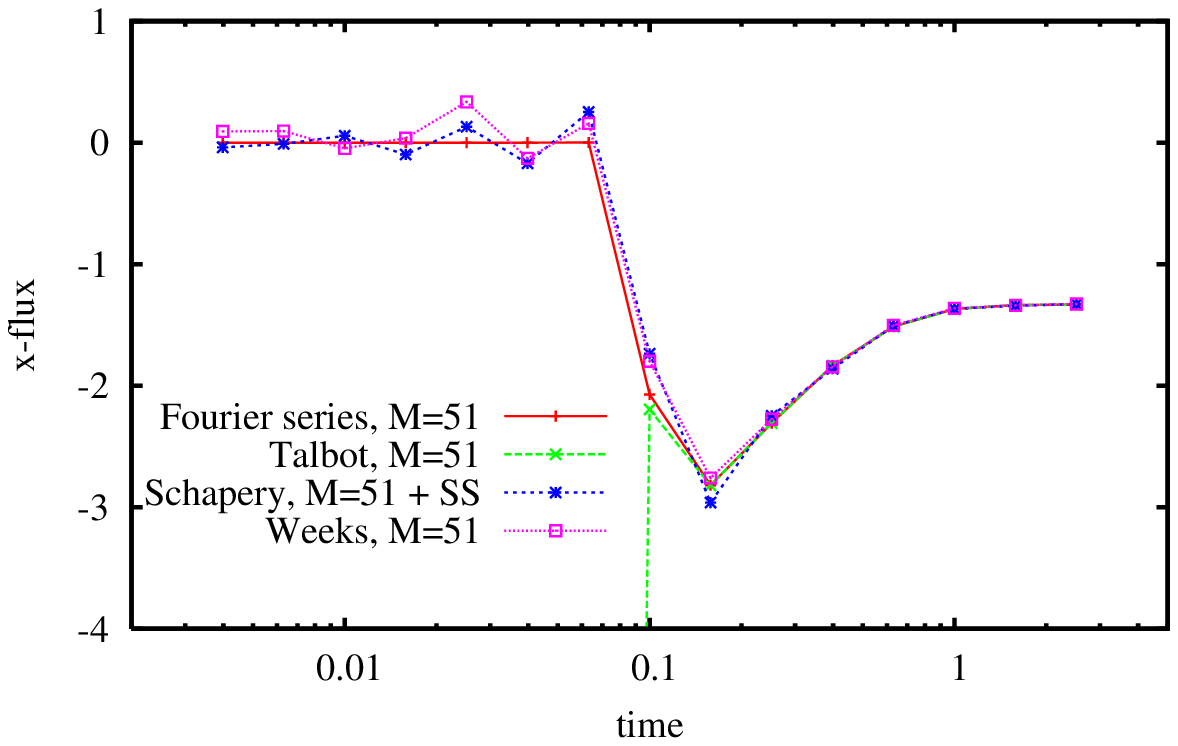}
  \caption{Plots of potential and flux through time with four methods
    for $\bar{f}_\mathrm{t}(p)=\exp(-0.08 p)/p$, same $p$
    used across all $t$. 51 total $\bar{f}(p)$ evaluations are used by each
    method.  Weeks' solution is undefined for $t<0.08$.}
  \label{fig:step-better}
\end{figure}

Finally, the same domain was simulated but with step-change boundary
conditions at $\tau=0.08$, or
$\bar{f}_\mathrm{t}(p)=\mathcal{L}(H(t-0.08))=e^{-0.08 p}/p$, with
singularities at the origin and $p=-\infty$. This function, and those
derived from it (e.g., a pulse or a square wave) are difficult
functions to invert accurately, because $f(t)$ is discontinuous.
Figures~\ref{fig:step-low} and \ref{fig:step-better} show the Talbot
method does not work for $t<\tau$ in double precision, since
$\bar{f}_\mathrm{t}(p)$ grows exponentially as $p \rightarrow
-\infty$. The Weeks and Schapery methods do worse than the
Fourier series approach (even with $N=51$), but their parameters can be
optimized further to improve these methods.

Although this step boundary condition could be implemented more
accurately by shifting the results from the first example by $t=0.08$,
other step-derived time behaviors involving a pulse or square wave
cannot be simplified in this way.

\subsection{Numerical Results Summary}
Table~\ref{tab:numericSummary} summarizes results from numerical
testing with these four simple boundary condition time behaviors.  The
second column indicates what limit there is on the number of terms in
the approximation and therefore the accuracy of the method. The size
of $p$ required by the Weeks and Fourier series methods grow much
slower than those required by the fixed Talbot method.  The third
column indicates what parameters are needed to be tuned by the
implementer to increase convergence of the method, and whether a good
choice is critical to the success of the method -- an automatic method
should not require searching or optimizing parameters to obtain a
robust solution.  We define robustness as the ability of a solution to
remain useful, even when not at optimality.  We prioritize a solution
that is good enough and stable over one that is excellent but
catastrophically sensitive to parameter choice.  The fourth column
indicates the ease of implementation in modern Fortran, Matlab, or
NumPy. The methods could also be implemented in a variable-precision
environment like Mathematica or mpmath \cite{mpmath}, but this would
further require the BEM model be implemented in such an environment.

\begin{table} 
  \caption{Numerical summary}
  \centering
  \begin{tabular}{|c|c|c|c|}
    \hline
    \textbf{Method} & \textbf{Number of Terms} &
    \textbf{Free Parameters}  & \textbf{Implementation} \\
    \hline
    Stehfest & $N\le$ decimal precision  & none & easiest \\
    Schapery & depends on choice of $p_j$ & $p_j$ via trial \& error & moderate \\
    Weeks & $p\rightarrow  i \infty$ slowly as $N$ grows  &
    $\kappa$ \& $b$ (very sensitive to $b$) & moderate \\
    fixed Talbot & $p \rightarrow -\infty$ quickly as $N$ grows &
    $r=\frac{2M}{5t_\mathrm{max}}$  (automatic) & easy \\
    Fourier series & $p \rightarrow i \infty$ slowly as $N$ grows  &
    $T=2t_\mathrm{max}$ (automatic) & most difficult \\
    \hline
  \end{tabular}
  \label{tab:numericSummary}
\end{table}

The modest success of the Schapery method is a bit surprising, given
its simplicity and use of real $p$. The results of the previous
section were the product of many iterations of trial and error, this
effort was not included in the implementation effort.  A better rule
or parametrization of $p_j$ might make this method more widely
useful.

The sensitivity of Weeks' method to the parameter choices was also
surprising; similarly, the method could have been improved after some
optimization \cite{weideman99}, but Weeks' rule of thumb was used for
the parameters. One of the noted advantages of Weeks' method is
the need to only compute optimal $p$ once, then any time can
accurately be inverted \cite{kano05,weideman99,duffy93}.  When using
the simple rules-of-thumb for the the free parameters, this was not
found to be the case.  The generalized form of Weeks' method can
include information about behavior of $\bar{f}(p) \rightarrow \infty$
(related to behavior as $t \rightarrow 0$), but this requires
problem-specific knowledge.

The Fourier series method is more robust with
respect to non-optimal $p$ values, even though \cite{duffy93}
cites this as a reason to use Weeks' method over the Fourier series
approach.

\section{Conclusions}
Laplace-space numerical approaches to solve the
diffusion equation have several viable alternative inverse Laplace
transform algorithms to choose from.  Historically, most Laplace-space
solutions to the diffusion equation have used real-only methods
(i.e., Gaver-Stehfest or Schapery). More robust methods require
complex arithmetic and $\bar{f}(p)$ evaluations, but have the benefits
of:
\begin{enumerate}
\item handling a broader class of time behaviors (Fourier series method);
\item still being relatively simple to implement (fixed Talbot method);
\item only utilizing double-precision complex data types, which are
  handled natively by modern Fortran, Matlab, or NumPy, and by common
  extensions in C++ (Fourier series and Weeks' methods).
\end{enumerate}

Several practical recommendations are made regarding Laplace-space
numerical modeling:
\begin{enumerate}
\item If many observations are needed across several time log cycles,
  large gains in efficiency can come from inverting groups of times
  with a single $\bar{f}(p)$ vector (e.g., grouped by log cycle).
  This complicates the implementation, but leads to much faster
  simulations.
\item If calculating $\bar{f}(p)$ is very expensive, and some
  numerical dispersion is allowable (not solving a wave problem with
  sharp fronts), then the Fourier series method approach is most
  economical, and is automatic and robust regarding free-parameter
  selection.
\item If only a single $\bar{f}_\mathrm{t}(p)$ is needed, then it may
  be worthwhile to optimize free parameters needed by Weeks' or
  Piessen's methods, or incorporate information about asymptotic
  $\bar{f}(p)$ behavior. Selection of optimum $b$ is far from
  automatic, and the Weeks method is not robust for non-optimal free
  parameters values.
\item If implementation time is a large factor, the fixed Talbot is
  quite simple to code and was automatic (no need to select optimum
  parameters).  The fixed Talbot may not work for non-zero step-time
  behavior without extended precision.
\item If complex-valued function evaluations are not feasible (e.g.,
  only real matrix or special function libraries are available), the
  Schapery or Piessen's methods are capable of using the same $p$
  values to invert different times, which the Gaver-Stehfest method
  cannot.
\item When appropriate, the strategy used by Schapery to expand the
  deviation from a reference state could be incorporated as a strategy
  to improve other algorithms.
\end{enumerate}
\end{linenumbers}

\begin{acknowledgements}
The author thanks Professors Cho Lik Chan and Barry Ganapol from the
Aerospace and Mechanical Engineering Department at the University of
Arizona for initial inspiration and direction on this manuscript.

Sandia National Laboratories is a multi-program laboratory managed and
operated by Sandia Corporation, a wholly owned subsidiary of Lockheed Martin
Corporation, for the U.S. Department of Energy's National Nuclear Security
Administration under contract DE-AC04-94AL85000.
\end{acknowledgements}


\end{document}